\documentclass[12pt]{amsart}
\usepackage[latin1]{inputenc}

\usepackage{xspace}
\usepackage{amsfonts,amssymb}
\usepackage{amsmath}
\usepackage{amsthm}
\usepackage{tikz}
\usepackage[latin1]{inputenc}
\usepackage[cyr]{aeguill}
\usepackage[english]{babel}
\usepackage{graphics}
\usepackage{vmargin}
\usepackage{xspace}
\usepackage{amsmath,amsfonts}

\def\footnoterule{\kern -15pt \hrule width 2truein \kern10.6pt\relax}%

\newtheorem{thm}{Theorem}[section]

\newtheorem{claim}{Claim}

\newtheorem{prop}[thm]{Proposition}
\newtheorem{df}[thm]{Definition}
\newtheorem{lm}[thm]{Lemma}
\newtheorem{cor}[thm]{Corollary}
\newtheorem{rmq}[thm]{Remark}
\newtheorem{prop-def}[thm]{Proposition and Definition}

\def\cercle{\mathbb{S}^1}
\def\H3{\mathbb{H}^3}
\def\zz{\mathbb{Z}}
\def\z2z{\mathbb{Z}/_{2\mathbb{Z}}}

\def\nn{\mathbb{N}}

\def\vol{{\rm Vol}}

\def\inj{{\rm Inj}}
\def\dist{{\rm dist}}
\def\d{{\rm d}}
\def\diam{{\rm diam}}
\def\lgr{{\rm length}\,}

\def\div{\rightarrow+\infty}

\newcommand{\incl}[1][r]
  {\ar@<-0.2pc>@{^(-}[#1] \ar@<+0.2pc>@{-}[#1]}

\title{Circular characteristics and fibrations of
 hyperbolic closed 3-manifolds.}

\author{Claire \textsc{Renard}.}

\begin{document}

\begin{abstract}

This article provides sufficient conditions for a closed
hyperbolic 3-manifold $M$ with non zero first Betti number to
fiber over the circle, and to find a fiber in $M$. Those
conditions are formulated in terms of the behavior the circular
characteristic in finite regular covers of $M$. We define the
circular characteristic as an invariant associated to a non
trivial cohomology class $\alpha$ of $M$, using a Heegaard
characteristic.

\end{abstract}

\maketitle

\today

\section*{Introduction}

Thurston conjectured that every complete hyperbolic, connected and
orientable 3-manifold of finite volume virtually fibers over the
circle, i.e. such a manifold has a finite covering that is a
bundle over the circle.

With this conjecture in mind, an interesting question is to find
criteria that are sufficient conditions for a closed hyperbolic
3-manifold $M$ to fiber over the circle. A necessary condition for
$M$ to be fibered is that its first Betti number $b_1(M)$ is
non zero.

The main idea of this article is to start with a non trivial
cohomology class $\alpha$ in $H^1(M,\zz)$ and to study the
behavior of a number associated to $\alpha$ called the circular
characteristic. This is a kind of Heegaard characteristic,
associated to a given non trivial cohomology class.

\begin{df}\label{def-deccirc}

Let $M$ be a hyperbolic, connected, oriented and closed
3-manifold. If $\alpha\in H^1(M)=H^1(M,\zz)$ is a non-trivial
cohomology class, let us denote by $\|\alpha \|$ the Thurston norm
of $\alpha$. By definition,
$$\|\alpha\|=\min\{\chi_-(R),\, [R]=\mathcal{P}(\alpha)\},$$
where $R$ is an embedded surface and $\mathcal{P}(\alpha)$ the
Poincar\'e-dual class of $\alpha$. We will call such a surface $R$
realizing the Thurston norm of $\alpha$ a
\textbf{$\|\alpha\|$-minimizing surface.}

If $R$ is a non-separating and $\|\alpha\|$-minimizing surface for
a given non-trivial cohomology class $\alpha\in H^1(M)$, take
$\mathcal{N}(R)\cong R\times (-1,1)$ a regular neighborhood of $R$
in $M$, and denote by $M_R=M\setminus \mathcal{N}(R)$. Set
$$h(M,\alpha,R)=\min\{\chi(R)-\chi(S)\},$$ where $S$ is a Heegaard
surface for $(M_R,R\times\{1\},R\times\{-1\})$. Said differently,
$\frac{1}{2}h(M,\alpha,R)$ is the minimal number of $1$-handles we
need to attach to a regular neighborhood of $R\times\{1\}$ in
$M_R$ to get the first compression body of a Heegaard splitting of
$(M_R,R\times\{1\},R\times\{-1\})$. Set $$h(\alpha)=h(M,\alpha)=
\min\{h(M,\alpha,R),\,[R]=\mathcal{P}(\alpha),\,
\chi_-(R)=\|\alpha\|\}.$$ For each non-trivial cohomology class
$\alpha\in H^1(M)$, let $\chi_-^c(\alpha)=\|\alpha\|+h(\alpha)$ be
the \textbf{circular characteristic} of $\alpha$. It is the
negative part of the Euler characteristic of a minimal genus
Heegaard surface for $M_R$, where $R$ is a $\|\alpha\|$-minimizing
surface such that the number $h(M,\alpha,R)$ is minimal among all
$\|\alpha\|$-minimizing surfaces.

\end{df}

The number $h(\alpha)$ can also be viewed as the minimal number of critical
points of a circular Morse function for $M$ such that the regular level sets
correspond to a surface the homology class of which is Poincar\'e dual to
$\alpha$. See section \ref{par-deccirc}.

\begin{rmq}\label{rmq-SH&DecCirc}

If $\alpha$ and $R$ are as above and $S$ is a Heegaard surface
corresponding to a Heegaard splitting of
$(M_R,R\times\{1\},R\times\{-1\})$ such that $\chi_-(R)=\| \alpha
\|$ and $\chi_-(S)=\chi_-^c(\alpha)$, then from the Heegaard
decomposition of $(M_R,R\times\{1\},R\times\{-1\})$, one can easily
construct a Heegaard splitting of $M$ by adding two small tubes
connecting the surfaces $R$ and $S$, each in one of the
compression bodies of the decomposition of
$(M_R,R\times\{1\},R\times\{-1\})$. An easy calculation shows that
\begin{align*}
\chi_-^h(M) &\leq \chi_-^c(\alpha) + \| \alpha \| +2\\
& \leq 2\chi_-^c(\alpha) +2.
\end{align*}

\end{rmq}

The idea is to use this number $\chi_-^c(\alpha)$ associated to a
given cohomology class $\alpha$ to get an explicit statement.
Studying the behavior of this circular characteristic when the
class $\alpha$ lifts to finite regular covers of $M$, we adapted
results of Lackenby \cite{Lac1} to obtain the following theorem,
which is the main result of this article.

\begin{thm}\label{thm-deccircregulier&fibration}

Let $M$ be a connected, oriented and closed hyperbolic
$3$-manifold, and set $\epsilon = \inj(M)/2$, where $\inj(M)$ is
the injectivity radius of $M$. 

There exists an explicit constant $\ell=\ell(\epsilon, \vol(M))$, depending only
on $\epsilon$ and the volume of the manifold $M$, and satisfying the following
properties. 

Let $\alpha\in H^1(M)$ be a non trivial cohomology
class and $R$ a $\|\alpha\|$-minimizing surface. Let
$M'\rightarrow M$ be a \textbf{regular} finite cover of $M$ of
degree $d$. Let $R'$ be a component of the preimage of $R$ in the
cover $M'$, and $\alpha'$ the cohomology class in $H^1(M',\zz)$
that is Poincar\'e-dual to $[R']$.

If $\ell\,\chi_-^c(\alpha') \leq \sqrt[4]{d},$ then the manifold
$M$ fibers over the circle and the surface $R$ is a fiber.

Furthermore, with $a'=6\left(\frac{21}{4}+\frac{3}{4\pi} +
\frac{3}{4\epsilon} +\frac{2}{\sinh^2(\frac{\epsilon}{4})}\right)$
and $D:=\frac{8\epsilon \vol(M)}{\pi(\sinh
(2\epsilon)-2\epsilon)}$, one has
$$\ell:=\sqrt[4]{117}\sqrt{a'\frac{\pi(\sinh(2D+2\epsilon)-2D-2\epsilon)}{2\vol(M)}}.
$$

\end{thm}

\begin{rmq}

The explicit expression of the constant $\ell$ involved in theorem
\ref{thm-deccircregulier&fibration} allows us to study its
behavior. If the volume $\vol(M)$ is fixed and that $\inj(M)$
tends to zero, or if $\inj(M)$ is fixed and $\vol(M)$ tends to
infinity, $\ell$ tends to infinity. Thus, the sufficient condition
given by the previous theorem becomes more and more difficult to
satisfy when the injectivity radius decreases (which corresponds
for example to a cusp opening), or if the volume grows (for
instance if one passes to finite covers of $M$).

\end{rmq}

The next corollary directly follows from theorem
\ref{thm-deccircregulier&fibration}.

\begin{cor}\label{cor-deccircregulier&fibration-asymptotique}

Let $M$ be a connected, oriented and closed hyperbolic
$3$-manifold. Let $\alpha\in H^1(M)$ be a non trivial cohomology
class and $R$ a $\|\alpha\|$-minimizing surface. Let
$(M_i\rightarrow M)_{i\in\nn}$ be a collection of finite regular
covers of $M$ with degrees $d_i$. For each $i\in\nn$, let $R_i$ be
a component of the preimage of $R$ in $M_i$, and $\alpha_i \in
H^1(M_i)$ the class that is Poincar\'e-dual to the class of $R_i$
in $H_2(M_i)$. If
$$\lim_{i\div}\frac{\chi_-^c(\alpha_i)}{\sqrt[4]{d_i}}=0,$$
then the manifold $M$ fibers over the circle, and the surface $R$
is a fiber.

\end{cor}

This corollary is true for any infinite collection of finite
covers satisfying the given asymptotic condition. If one considers
the tower of cyclic finite covers of $M$ dual to the class
$\alpha$, theorem \ref{thm-deccircregulier&fibration} leads to the following
corollary.

\begin{cor}\label{cor-revcycliquesduauxetcl}

let $M$ be a connected, oriented and closed hyperbolic
$3$-manifold. Let $\alpha\in H^1(M)$ be a non trivial cohomology
class and $R$ a $\|\alpha\|$-minimizing surface. Let
$(M_i\rightarrow M)_{i\in\nn}$ be the collection of cyclic finite
covers of $M$ dual to the class $\alpha$, such that for every
$i\in\nn$, the cover $p_i\,:\,M_i\rightarrow M$ is regular, with
degree $i$. For each $i\in\nn$, let $\alpha_i:=p_i^*(\alpha)$ be
the cohomology class in $H^1(M_i,\zz)$ corresponding to $\alpha$.

If there exists $i\geq i_0=\lceil (2\ell\|\alpha\|)^4\rceil$ such that
$$\frac{h(\alpha_i)}{\sqrt[4]{i}}\leq \frac{1}{4\ell},$$
then the manifold $M$ fibers over the circle, and the surface $R$
is a fiber .

\end{cor}

~

\noindent \textbf{Acknowledgement:} I would like to thank warmly my advisor,
Michel Boileau, whose encouragements, kindness and patience were
essential ingredients in this work. I am grateful to Juan Souto,
Nicolas Bergeron, Joan Porti, Jean-Marc Schlenker, Jean-Pierre Otal, Vincent
Guirardel and Cyril Lecuire for very helpful conversations during the
elaboration of this paper.

\section{Circular decompositions and thin position.}\label{par-deccirc}

A circular decomposition is the equivalent of a Heegaard
decomposition, but this decomposition is associated to a Morse
function that no longer takes values in $I=[0,1]$ but in the
circle $\cercle$. According to \cite{MG}, we have the following definitions.

\begin{df}

A \textbf{circular Morse function} is a Morse function $f\,:\, M
\rightarrow \cercle$.

If $f\,:\, M \rightarrow \cercle$ is a circular Morse function,
the handle decomposition of $M$ given by the function $f$ is
called the \textbf{circular decomposition associated to $f$}.

\end{df}

See F. Manjarrez-Guti\'errez \cite{MG}, Matsumoto \cite{Mat} and
Milnor \cite{Mil} for further details about circular Morse
functions. Let $f\,:\,M\rightarrow \cercle$ be a circular Morse
function. If we remove a small open neighborhood of a regular
value $x \in \cercle$, by restriction of $f$, we obtain a Morse
function $g$ of $M_R=M\setminus \mathcal{N}(R)$, which is the
manifold $M$ minus a small regular open neighborhood of the
surface $R:=f^{-1}(\{x\})$, on the interval $I$. Thus, the theory
of Heegaard splittings and generalized Heegaard splittings applies
to the function $g$.

~

An other viewpoint is to see a circular decomposition as a handle
decomposition of the cobordism $(M\setminus\mathcal{N}(R),
R\times\{1\}, R\times \{-1\})$. Starting with a Heegaard splitting
of Heegaard surface $S$ for $M_R=M\setminus\mathcal{N}(R)$, one
can change the order in which $1$- and $2$-handles are attached to
get a new generalized Heegaard splitting
$(F_1=R\times\{1\},S_1,F_2,\ldots,S_n,F_{n+1}=R\times\{-1\})$ for
$(M_R, R\times\{1\}, R\times \{-1\})$. Gluing back $R\times\{1\}$
to $R\times\{-1\}$, one obtains a circular decomposition for the
manifold $M$. Denote it by
$\mathcal{H}=(F_1,S_1,F_2,\ldots,S_n,F_{n+1})$, with
$F_1=F_{n+1}=R$. The surfaces $F_j$ divide $M$ into $n$
3-manifolds with boundary $W_1,\ldots,W_n$, and surfaces $S_j$ are
Heegaard surfaces for those manifolds. For $1\leq j\leq n$, $S_j$
divides the manifold $W_j$ into two compression bodies $A_j$ and
$B_j$, such that $\partial_+A_j=\partial_+B_j=S_j$,
$\partial_-A_j=F_j$ and $\partial_-B_j=F_{j+1}$.

Let $S$ be a closed surface. If $S$ is connected, recall that the
\textbf{complexity} of $S$ is $c(S)=\max(0,2g(S)-1)$. If $S$ is
the union of several connected components, the complexity of $S$
is the sum of the complexities of the components of $S$. There is
a definition of the complexity of a circular decomposition
analogous to the complexity of a generalized Heegaard splitting.

\begin{df}

The \textbf{circular width} of a circular decomposition
$\mathcal{H}=(F_1,S_1,F_2,\\ \ldots,S_n,F_{n+1})$ is the set of
the $n$ integers $(c(S_1),\ldots,c(S_n))$, with repetitions and
arranged in monotonically non-increasing order. Widths are
compared using the lexicographic order.

The integer $n\geq 1$ is called the \textbf{length} of the
circular decomposition
$\mathcal{H}=(F_1,S_1,\\F_2,\ldots,S_n,F_{n+1})$.

\end{df}

\begin{prop}\label{prop-dec-mince}

Let $M$ be a hyperbolic, connected, oriented and closed
3-manifold. Let $R$ be an orientable, closed, non-separating,
incompressible and embedded surface in $M$. Denote by $S$ a
Heegaard surface for $M\setminus \mathcal{N}(R)$. Starting from
the circular decomposition $\mathcal{H}=(R,S,R)$ of $M$, there
exists a finite number of surgeries to get a circular
decomposition $\mathcal{H'}=(F_1,S_1,F_2,\ldots,S_n,F_{n+1})$ with
$F_1=F_{n+1}=R$, such that:
\begin{enumerate}
  \item the circular width of $\mathcal{H'}$ is minimal among the
  widths of such circular decompositions obtained by a finite number
  of surgeries of $\mathcal{H}$,
  \item each surface $S_j$ is a strongly irreducible Heegaard
  surface for the Heegaard decomposition $(A_j,B_j)$ of $W_j$
  and $g(S_j)\leq g(S)$,
  \item each surface $F_j$ is incompressible, no component of $F_j$
  is a sphere, and $g(F_j)\leq g(S)$,
  \item $n\leq \frac{1}{2}(\chi(R) - \chi(S))$,
  \item $\chi(R)-\chi(S)=\sum_{j=1}^n( \chi(F_j) -\chi(S_j))$.
  \item Furthermore, if the decomposition $\mathcal{H'}$ is of
  length at least $2$, up to forgetting some surfaces, one can
  assume that for every $j$, the surfaces $F_j$ and $F_{j+1}$ are
  not parallel.
\end{enumerate}

\end{prop}

\begin{df}

Let $\mathcal{H}$ be a circular decomposition. A circular
decomposition $\mathcal{H'}=(F_1,S_1,F_2,\ldots,S_n,F_{n+1})$ that
is circular-length-minimizing among all circular decompositions
obtained from $\mathcal{H}$ by a finite number of surgeries is
said to be a \textbf{thin position}. We will call such a
decomposition a \textbf{thin circular decomposition associated to
$\mathcal{H}$.}

\end{df}

\noindent \underline{Proof of proposition \ref{prop-dec-mince}.}

The proof of the first three points of this proposition is based on the proof of
\cite[Theorem 3.2]{MG}, which is itself an adaptation of
techniques of \cite{ST2} to the case of circular decompositions. We recall here
the arguments (see also \cite[section 3]{Lac}).

We start with the circular decomposition $\mathcal{H}=(R,S,R)$.
The aim is to perform a certain number of surgeries to obtain a
decomposition in a thin position, i.e. of minimal complexity. Each
surgery corresponds to a change on the order in which $1$- and
$2$-handles are attached, such that a surgery strictly decreases
the circular width of the decomposition. Thus, the number of
necessary surgeries to get a thin decomposition is finite.

\begin{lm}\label{lm-SHfaiblt-red->chirurgie}

Let $\mathcal{H}=(F_1,S_1,\ldots,S_n,F_{n+1})$ be a circular
decomposition for $M$, and suppose that for some index $j$, the
Heegaard surface $S_j$ for $(A_j,B_j)$ is weakly reducible. Then
there exists an operation called a \textbf{surgery}, starting from
$\mathcal{H}$ and giving a new circular decomposition
$\mathcal{H'}$ of strictly smaller circular width.

\end{lm}

\noindent \underline{Proof of lemma
\ref{lm-SHfaiblt-red->chirurgie}.}
\nopagebreak

As the Heegaard surface $S_j$ is weakly reducible, there exists a
pair of disjoint compression discs for $S_j$, say $D_A$ embedded
in $A_j$ and $D_B$ in $B_j$. Performing surgeries along those two
discs, one gets a new circular decomposition $\mathcal{H'} :=
(F_1,\ldots,F_j, T_j,G_j,\\T_j',F_{j+1},\ldots,F_{n+1})$, where
the surface $T_j$ is obtained from $S_j$ by surgery along $D_A$,
$T_j'$ from $S_j$ by surgery along $D_B$, and $G_j$ from $S_j$ by
surgery along $D_A$ and $D_B$. As $\left|\chi(T_j) \right| =
\left|\chi(T'_j) \right| =\left|\chi(S_j) \right| -2$, the
circular width of this new circular decomposition is strictly
smaller than this of $\mathcal{H}$.

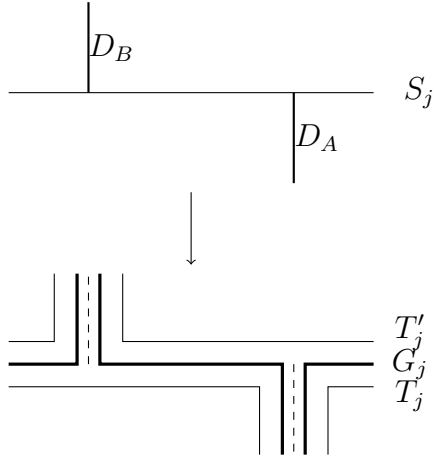
\begin{center}
\begin{tikzpicture}[scale=.6]

\draw (-4,6) -- (4,6);

\draw [thick] (-2.25,6) -- (-2.25,8);

\draw [thick] (2.25,6) -- (2.25,4);

\node at (5,6) {$S_j$};

\node at (-1.75,7) {$D_B$};

\node at (2.75,5) {$D_A$};

\draw [->] (0,3.8) -- (0,2.2);

\draw [very thick] (-4,0) -- (-2.5,0) -- (-2.5,2);

\draw [very thick] (-2,2) -- (-2,0) -- (2,0) -- (2,-2);

\draw [very thick] (2.5,-2) -- (2.5,0) -- (4,0);

\draw (-4,.5) -- (-3,.5) -- (-3,2);

\draw (-1.5,2) -- (-1.5,.5) -- (4,.5);

\draw (4,-.5) -- (3,-.5) -- (3,-2);

\draw (1.5,-2) -- (1.5,-.5) -- (-4,-.5);

\draw [dashed] (-2.25,0) -- (-2.25,2);

\draw [dashed] (2.25,0) -- (2.25,-2);

\node at (4.8,0) {$G_j$};

\node at (4.8,.7) {$T_j'$};

\node at (4.8,-.7) {$T_j$};

\end{tikzpicture}
\end{center}\qed

~

As $\chi(T_j)=\chi(T_j')=\chi(S_j)+2$ and $\chi(G_j)=\chi(S_j)+4$,
one obtains $-\chi(S_j)=-\chi(T_j)+\chi(G_j)-\chi(T_j')$. Thus,
this surgery procedure does not modify the alternate sum
$\sum(\chi(F_j)-\chi(S_j))$, proving point $(5)$.

As this surgery procedure strictly decreases the circular width of
the decomposition, there exists a finite number of such surgeries
to get a circular decomposition
$\mathcal{H'}=(F_1,S_1,\ldots,S_n,F_{n+1})$ of minimal circular
width among the set of all decompositions obtained by surgeries
from the starting circular decomposition $\mathcal{H}=(R,S,R)$.

To prove $(2)$, recall \cite{MG}. Just notice that if one of the Heegaard
surfaces $S_j$ is not strongly irreducible, from lemma
\ref{lm-SHfaiblt-red->chirurgie}, one can perform another surgery
to obtain a new circular decomposition of circular width strictly
smaller than this of $\mathcal{H'}$, which is a contradiction if
$\mathcal{H'}$ is a length-minimizing decomposition.

The proof of point $(3)$ is done in \cite{MG}. The
surface $R=F_1=F_{n+1}$ is incompressible. Suppose by
contradiction that one of the surfaces $F_j$ is compressible, for
an index $j$ between $2$ and $n$. There exists then a compression
disc $D$ for $F_j$. Taking an innermost disc, one can furthermore
assume that $D\cap (\cup_{k=1}^n F_k)=D\cap F_j=\partial D$. Thus,
the disc $D$ entirely lies in the region $W_{j-1}$ bounded by the
two surfaces $F_{j-1}$ and $F_j$, or is entirely embedded in the
region $W_j$ bounded by $F_j$ and $F_{j+1}$. Suppose for example
that $D$ is entirely embedded in $W_j$. From the boundary version
of the Haken Lemma \cite{Hak}, as $W_j$ is $\partial$-reducible,
every Heegaard splitting of $W_j$ is reducible, hence weakly
reducible. This is a contradiction with point $(2)$, proving the
first part of point $(3)$.

If one of the components of a surface $F_j$ is a 2-sphere, as $M$
is irreducible, this sphere bounds an embedded ball in $M$. Taking
an innermost sphere, one obtains a sphere bounding the Heegaard
splitting of a 3-ball. But this splitting, if not trivial, is
reducible (see \cite{Wald}), hence weakly reducible, contradicting point $(2)$.
This ends the proof of point $(3)$.

To prove point $(4)$, notice that the surgery process as described
above is in fact a change on the order in which the handles are
attached. More precisely, with the notations above, if we consider
a handle decomposition associated to $\mathcal{H}$ where $1$- and
$2$-handles correspond to meridian discs for the Heegaard
splittings, a surgery is a handle reordering. The $2$-handle
corresponding to the meridian disc $D_B$ is attached before the
$1$-handle corresponding to $D_A$. Thus, this process does not
change the number of $1$- and $2$-handles. In the starting
circular decomposition $\mathcal{H}=(R,S,R)$, the number of $1$-
and $2$-handles is equal to $\chi(R) - \chi(S)$. So after each
surgery, there are still $\frac{1}{2}(\chi(R) - \chi(S))$
$1$-handles and $\frac{1}{2}(\chi(R) - \chi(S))$ $2$-handles. As
the number of regions of a circular decomposition $\mathcal{H'}$
is at most the number of $1$- and $2$-handles in this
decomposition, there are at most $\chi(R) - \chi(S)$ regions in
$\mathcal{H}$. Therefore, the number of even surfaces $F_j$ is
bounded above by $\frac{1}{2}(\chi(R) - \chi(S))$. In other words,
$n\leq \frac{1}{2}(\chi(R) - \chi(S))$, which proves point $(4)$.

Eventually, for point $(6)$ we recall the argument of \cite[Section 3]{Lac}. If
the length of the decomposition is just $1$, this means that there is only one
incompressible surface $F_1=R=F_2$. If $F_1$ is parallel to $F_2$ in $M_R$, in
fact the manifold $M$ fibers over the circle, with fiber $R$.

If the length of the decomposition $\mathcal{H'}$ is at least $2$,
suppose that there exists two parallel surfaces $F_j$ and
$F_{j+1}$ for some $j$. From point $(2)$, the surface $S_j$ is a
strongly irreducible Heegaard surface for the product region
bounded by $F_j$ and $F_{j+1}$. From the classification of
Heegaard splittings for products (see \cite{ST1}), this means that
$S_j$ is parallel to $F_j$. The two surfaces $F_j$ and $F_{j+1}$
can then be amalgamated to a single surface, forgetting the
surface $S_j$, to obtain a new circular decomposition with
complexity strictly smaller than this of $\mathcal{H'}$ and still
verifying the other points of proposition
\ref{prop-dec-mince}.\qed

\begin{cor}\label{cor-estimates}

Let $M$, $R$ and $S$ be as above, and
$\mathcal{H'}=(F_1=R,S_1,\ldots, F_{n+1}=R)$ a thin circular
decomposition associated to $(R,S,R)$. Let $\overline{F}$ be the
surface obtained from $\bigcup_jF_j\cup \bigcup_jS_j$ by
amalgamating parallel components bounding product regions in $M
\setminus \left( \bigcup_jF_j\cup \bigcup_jS_j \right)$ in a
single component. Then,
\begin{enumerate}
  \item $\left|\chi(\overline{F})\right|\leq \left|\chi(\bigcup_jF_j\cup
  \bigcup_jS_j)\right|\leq \left|\chi(S)-\chi(R)\right| \left| \chi(S)
  \right|$,
  and
  \item the surface $\overline{F}$ has at most
  $\frac{3}{2}\left|\chi(S)-\chi(R)\right|$ connected components.
\end{enumerate}

\end{cor}

\noindent \underline{Proof of corollary \ref{cor-estimates}.}
\nopagebreak

We adapt here the proof of \cite[Corollary 4]{Lac1}. First, notice
that no compression body in the complement of $\bigcup_jF_j\cup
\bigcup_jS_j$ is a punctured $3$-ball, as no
component of $\bigcup_jF_j\cup \bigcup_jS_j$ is a $2$-sphere.

As $M$ is hyperbolic, no compression body of the thin circular decomposition can
be a solid torus.

\begin{rmq}

An other way to prove point $(4)$ of proposition \ref{prop-dec-mince} starting
from point $(5)$ is the following.

Recall that $F_1=R=F_{n+1}$. Point $(5)$ of proposition
\ref{prop-dec-mince} can also be written:
\begin{equation}\label{eq-chi}
\chi(R)-\chi(S)=\frac{\chi(F_1)-\chi(S_1)}
{2}+\frac{\chi(F_2)-\chi(S_1)}{2} +\frac{\chi(F_2)-\chi(S_2)}{2}+
\ldots +\frac{\chi(F_{n+1})-\chi(S_n)}{2}.
\end{equation}

If $H$ is a compression body that is not a punctured $3$-ball, nor
a solid torus, nor a product, then
$\chi(\partial_-H)-\chi(\partial_+H)>0$ and this integer is even.
As the $2n$ components of the complementary of $\bigcup_jF_j\cup
\bigcup_jS_j$ are such compression bodies, the right hand side of
equality (\ref{eq-chi}) is bounded from below by $2n$. Therefore,
$2n\leq \chi(R)-\chi(S)$. It is exactly point (4) of proposition
\ref{prop-dec-mince}.

\end{rmq}

Thus,
\begin{eqnarray*}
\left|\chi(\bigcup_jF_j\cup \bigcup_jS_j)\right| &=&
\sum_{j=1}^n \left|\chi(F_j)\right| + \sum_{j=1}^n \left|\chi(S_j)\right|\\
&\leq & 2n\left|\chi(S)\right|\\
&\leq & \left|\chi(R)-\chi(S)\right| \left|\chi(S)\right|.
\end{eqnarray*}

As some components of $\bigcup_jF_j\cup \bigcup_jS_j$ have been
discarded to form the surface $\overline{F}$, this implies
$\left|\chi(\overline{F}) \right|\leq \left| \chi(\bigcup_jF_j\cup
\bigcup_jS_j)\right|$, which proves point $(1)$ of corollary
\ref{cor-estimates}.

If $H$ is a compression body that is not a punctured $3$-ball, nor
a solid torus, nor a product, one can check that $\left|
\partial H \right|\leq
\frac{3}{2}(\chi(\partial_-H)-\chi(\partial_+H))$. The sum over
all compression bodies $H$ in the complement of $\bigcup_jF_j\cup
\bigcup_jS_j$ of $\chi(\partial_-H) - \chi(\partial_+H)$ is equal
to $\sum_H(\chi(\partial_-H)-\chi(\partial_+H)) = 2 \sum_{j=1}^n(
\chi(F_j) -\chi(S_j)) = 2(\chi(R)-\chi(S))$. Now, the number of
components of  $\overline{F}$ is at most $\frac{1}{2} \sum_{H}
\left|\partial H\right|$, where $H$ describes all compression
bodies that are the components of $M\setminus
\left(\bigcup_jF_j\cup \bigcup_jS_j\right)$ which are not product regions. But
\begin{eqnarray*}
\frac{1}{2} \sum_{H} \left|\partial H\right| &\leq & \frac{1}{2} \sum_H \frac{3}{2}(\chi(\partial_-H)-\chi(\partial_+H))\\
&=& \frac{3}{2} \sum_{j=1}^n( \chi(F_j) -\chi(S_j))\\
& = &\frac{3}{2}\left|\chi(R)-\chi(S)\right|.
\end{eqnarray*}
Therefore, $\left|\overline{F}\right|\leq
\frac{3}{2}\left|\chi(R)-\chi(S)\right|$, which ends the proof of
corollary \ref{cor-estimates}.\qed

~

The proof of theorem \ref{thm-deccircregulier&fibration} will
require to control the metric of the surface $\bigcup_jF_j\cup
\bigcup_jS_j$ of a thin circular decomposition of the hyperbolic
manifold $M$.

\begin{df}

An embedded surface $S$ in a Riemannian 3-manifold $M$ is called
\textbf{pseudo-minimal} if it is orientable, closed, and $S$ is a
minimal surface or the boundary of a regular neighborhood of a
minimal non-orientable surface, possibly with a little tube
attached vertically in the $I$-bundle structure.

\end{df}

Part (1) of the following theorem comes from results of Frohman,
Freedman, Hass and Scott about incompressible surfaces (\cite{FHS}
and \cite{FrHa}). Part (2) is a result of Pitts and Rubinstein
(\cite{PiRu}, see also \cite[Existence Theorem of minimal surfaces]{Sou3},
\cite{CoDL} and \cite{Pi}).

\begin{thm}\label{thm-surf-min}

Let $N$ be a connected, oriented and closed hyperbolic 3-manifold.

$(1)$ Any incompressible surface in $N$ can be isotoped to a
minimal surface or the boundary of a small neighborhood of a
non-orientable minimal surface.

$(2)$ Any embedded surface corresponding to a strongly irreducible
Heegaard surface for a region of $N$ lying between two (possibly
empty) embedded, incompressible and pseudo-minimal surfaces as
above can be isotoped to a minimal surface, or to the boundary of
a small regular neighborhood of a non-orientable minimal surface,
with a small tube attached vertically in the $I$-bundle structure.

\end{thm}\qed

~

The next corollary directly follows from theorem
\ref{thm-surf-min} combined with proposition \ref{prop-dec-mince}.

\begin{cor}\label{cor-dec-mince&surfPM}

Let $M$ be a hyperbolic, connected, oriented and closed
3-manifold. Take $\mathcal{H}=(F_1,S_1,F_2,\ldots,S_n,F_{n+1})$ a
thin circular decomposition of $M$. Then, up to isotopy, one can
assume that all surfaces $F_j$ and $S_j$ are pseudo-minimal.

\end{cor}\qed

\section{Homology classes and fibration of finite regular covers.}

The aim of this section is to prove theorem
\ref{thm-deccircregulier&fibration} and corollaries
\ref{cor-deccircregulier&fibration-asymptotique} and
\ref{cor-revcycliquesduauxetcl}.

\noindent \underline{Proof of theorem
\ref{thm-deccircregulier&fibration}.}
\nopagebreak

The proof is an adaptation of the proof of \cite[Theorem 1
(3)]{Lac1}, together with some calculations of explicit constants.

Let $M$ be a connected, oriented and closed hyperbolic 3-manifold
as in the assumptions of theorem
\ref{thm-deccircregulier&fibration}, and $\epsilon \leq
\inj(M)/2$. Let $\alpha\in H^1(M)$ be a non trivial cohomology
class and $R$ a $\|\alpha\|$-minimizing surface. Let
$M'\rightarrow M$ be a \textbf{regular} finite cover of $M$ with
degree $d$. Let $R'$ be a connected component of the preimage of
$R$ in the cover $M'$, and $\alpha'$ the cohomology class in
$H^1(M',\zz)$ that is Poincar\'e-dual to $[R']$. First, for needs
of the proof, suppose that in addition the surface $R'$ is
$\|\alpha'\|$-minimizing and such that
$h(M',\alpha',R')=h(M',\alpha')$.

Let $S'$ be a minimal genus Heegaard surface for
${M'}_{R'}=M'\setminus \mathcal{N}(R')$. By construction,
$\left|\chi(S')\right|=\chi_-^c(\alpha')$. Applying proposition
\ref{prop-dec-mince} to the circular decomposition $(R',S',R')$,
one obtains a thin circular decomposition
$\mathcal{H'}=(F_1,S_1,\ldots, F_{n_i+1})$ associated to
$(R',S',R')$, where $F_1=F_{n_i+1}=R'$. Moreover, all surfaces
$F_j$ and $S_j$ are isotopic to pseudo-minimal surfaces. If
$\overline{F}$ is the surface obtained from $\bigcup
F_j\cup\bigcup S_j$ as in corollary \ref{cor-estimates}, then
$\overline{F}$ is a pseudo-minimal surface, and it follows from
corollary \ref{cor-estimates} that
$\left|\chi(\overline{F})\right|\leq \left|
\chi(R')-\chi(S')\right| \left|\chi(S'\right|\leq
\chi_-^c(\alpha')^2$.

Let $\mathcal{D}$ be a Dirichlet fundamental polyhedron for the
manifold $M$ in its universal cover $\H3$. The union of the
translates of $\mathcal{D}$ under the action of the fundamental
group of $M$ composes a tiling of $\H3$. By the covering map $\H3
\rightarrow M'$, this tiling projects to a tiling of $M'$ by $d$
copies of $\mathcal{D}$. Let $\mathcal{D'}$ be one of those
polyhedra. As the cover $M' \rightarrow M$ is regular, the tiling
of $M'$ is the union of the translates of $\mathcal{D'}$ under the
action of the group $G:=\pi_1(M)/\pi_1(M')$.

One needs a few definitions and lemmas.

\begin{df}

Let $\epsilon>0$. The \textbf{$\epsilon$-diameter} of a
non-necessarily connected surface $F$ is the minimal number of
balls of radius $\epsilon$ for the metric of $F$ required to cover
the surface $F$.

\end{df}

\begin{lm}\label{lm-controle-diam-surf-min}

Let $S$ be an embedded pseudo-minimal surface in $N$, a Riemannian closed
3-manifold, whose sectional curvature is at most $-1$. Let
$\epsilon\leq \inj(N)/2$ and
$$a'=6\left(\frac{21}{4}+\frac{3}{4\pi}+ \frac{3}{4\epsilon}+
\frac{2}{\sinh^2(\frac{\epsilon}{4})} \right) .$$

Then the $\epsilon$-diameter of the surface $S$ is bounded from
above by $a'\left|\chi(S)\right|/3$. Furthermore, the
$\epsilon$-diameter of a pseudo-minimal surface $\Sigma$ homotopic
to $S$ and close enough is at most $a'\left|\chi(\Sigma)\right|$.

\end{lm}

\noindent \underline{Proof of lemma \ref{lm-controle-diam-surf-min}.}
\nopagebreak

This lemma is a direct consequence of \cite[Lemma 4.2 p.
2249]{Mah} and \cite[Proposition 6.1]{Lac} in the case the surface $S$ is
minimal and orientable, and we can take $a'/6$ instead of $a'$.
If $S$ is minimal, but not orientable, its homology class $[S]$ is non zero in
$H_2(N,\zz/2\zz)$. By Poincar\'e's duality, it corresponds to a
non-trivial element $\alpha\in H^1(N,\zz/2\zz)$. As the homology
class of the double cover of $S$ can be represented by the
boundary of a small regular neighborhood of the non-orientable
surface $S$, we have $2[S]=0$ in $H_2(N,\zz)$. If we take the
double cover $N'$ of $N$ corresponding to the kernel of $\alpha$,
the surface $S$ lifts to a minimal orientable surface $S'$. We can
apply the stronger version of lemma
\ref{lm-controle-diam-surf-min}, and bound the $\epsilon$-diameter
of $S'$ by $a'/6\left|\chi(S')\right| =a'/6
\times2\left|\chi(S)\right|= a'/3 \left|\chi(S)\right|$, and the
length of a one-vertex triangulation for $S'$ by $2\epsilon
a'/3\left|\chi(S)\right|$. As those numbers bound also from above
the $\epsilon$-diameter and the length of a one-vertex
triangulation of $S$, this proves the lemma for a minimal non orientable
surface, with $a'/3$ instead of $a'$.

If the surface $S$ is just pseudo minimal, it is the boundary of an arbitrarily
small regular neighborhood of a minimal surface $S'$. As the diameter and the
length of the edges of a one-vertex triangulation are at most $a'/3 \left|
\chi(S')\right|$ and $2\epsilon a'/3 \left| \chi(S') \right|$, with
$\left|\chi(S) \right| \leq 2 \left| \chi(S') \right|$, this ends the proof
of lemma \ref{lm-controle-diam-surf-min}. \qed

~

The following lemma is a way to
bound the diameter of a fundamental polyhedron $\mathcal{D}$ in
$\H3$ in terms of the volume of the manifold $M$ and a lower bound
for its injectivity radius.

\begin{lm}\label{lm-estimation_Dalpha&L}

Let $\mathcal{D}$ be a Dirichlet fundamental polyhedron for the
manifold $M$, embedded in the universal cover $\widetilde{M}\simeq
\mathbb{H}^3$. Let $D$ be an upper bound for the diameter of
$\mathcal{D}$ in $\H3$. We have the following estimate:

\begin{equation}\label{estimation_D}
\diam(\mathcal{D})\leq \frac{8\epsilon \vol(M)}{\pi(\sinh
(2\epsilon)-2\epsilon)}=D.
\end{equation}

If $S$ is an embedded surface in the finite cover $M'$ of $M$,
which can be covered by at most $\lambda$ embedded balls in $M'$
of radius $\epsilon\leq \inj(M)$, then $S$ intersects at most $L$
images of $\mathcal{D}$ in $M'$, with

\begin{equation}\label{estimation_L}
L=\lfloor \frac{\pi(\sinh(2D+2\epsilon)-2D-2\epsilon)}
{\vol(M)}\lambda\rfloor.
\end{equation}

\end{lm}\qed

\noindent \underline{Proof of lemma \ref{lm-estimation_Dalpha&L}.}
\nopagebreak

To prove inequality (\ref{estimation_D}), first notice that
$\diam(\mathcal{D})\leq 2\,\diam(M)$. To prove it, recall that
there exists $w\in\mathbb{H}^3$ such that
$\mathcal{D}=\{x\in\mathbb{H}^3\,,\,\d(\gamma(w),x)\geq\d(w,x)\;\forall
\gamma\in \pi_1(M)\}$. If $x$ and $y\in\mathcal{D}$ satisfy
$\d(x,y)=\diam(\mathcal{D})$, then
$$\diam(\mathcal{D})=\d(x,y)\leq\d(x,w)+\d(y,w)\leq 2\,\diam(M).$$

Take $x$ and $y\in M$ such that $\d(x,y)=\diam(M)$, and let
$\gamma$ be a minimizing geodesic from $x$ to $y$. By definition,
$\lgr(\gamma)=\diam(M)$. Let $\mathcal{B}$ be a collection of
points in $\gamma$ which is maximal among collections of points of
$\gamma$ such that two balls of radius $\epsilon$ and whose
centers are two distinct points of $\mathcal{B}$ have disjoint
interiors. Then, by maximality of $\mathcal{B}$, the union of
balls with centers in $\mathcal{B}$ and radius $2\epsilon$ cover
the geodesic $\gamma$.

Thus, $\left|\mathcal{B}\right| \geq \frac{\lgr(\gamma)}
{4\epsilon}$. As balls of centers in $\mathcal{B}$ and radius
$\epsilon$ have disjoint interiors, considering volumes, we
deduce:

\begin{eqnarray*}
\vol(M)&\geq & \sum_{u\in\mathcal{B}}\vol(B(u,\epsilon))\\
&\geq & \frac{\lgr(\gamma)}{4\epsilon} \vol(B_{\mathbb{H}^3}(\epsilon))\\
&\geq &\frac{\diam(M)}{4\epsilon} \pi(\sinh(2\epsilon)
-2\epsilon),
\end{eqnarray*}
proving inequality (\ref{estimation_D}).

To prove inequality (\ref{estimation_L}), denote by
$\mathcal{B}$ the set of the centers of a collection of $K$
embedded balls in $M'$ of radius $\epsilon'$ covering the surface
$S$. Let $\mathcal{N}=\cup_{x\in\mathcal{B}}B(x,D+\epsilon')$.
Those balls are not necessarily isometric to hyperbolic embedded
balls in $\H3$ as $D+\epsilon'>\inj(M)$. However, let us show that
$\mathcal{N}$ contains every fundamental polyhedron of $M'$
intersecting $S$.

To prove it, let $x$ be a point in a fundamental polyhedron of
$M'$ intersecting $S$. Take $y\in S$ such that
$\d(x,y)=\dist(x,S)\leq D$. As $y$ is a point of $S$, there exists
a ball $B(x,\epsilon')$ with $x\in \mathcal{B}$ containing $y$.
Therefore $\d(z,x)\leq \d(z,y)+\d(y,x) \leq D+\epsilon'$, showing
that $z\in B(x,\epsilon'+D)\subset \mathcal{N}$.

Comparing volumes, we get:
\begin{eqnarray*}
L\,\vol(\mathcal{D})&\leq &\vol(\mathcal{N})\\
L\,\vol(M)&\leq & \left|\mathcal{B}\right|\vol(B_{\H3}(\epsilon'+D))\\
L&\leq
&\frac{\pi(\sinh(2\epsilon'+2D)-2\epsilon'-2D)}{\vol(M)}K,
\end{eqnarray*}
proving inequality (\ref{estimation_L}), as
$L$ is a natural integer.\qed

~

In the sequel, set $a'=6\left(\frac{21}{4}+\frac{3}{4\pi} +
\frac{3}{4\epsilon} +\frac{2}{\sinh^2(\frac{\epsilon}{4})}\right)$
and $D:=\frac{8\epsilon \vol(M)}{\pi(\sinh
(2\epsilon)-2\epsilon)}$ as in lemmas
\ref{lm-controle-diam-surf-min} and \ref{lm-estimation_Dalpha&L}.
As $D$ is an upper bound for the diameter of $\mathcal{D}$ in
$\H3$, it is also an upper bound for the diameter of
$\mathcal{D'}$ in $M'$.

\begin{lm}\label{lm-intersectionD'}

Set
$\kappa:=a'\frac{\pi(\sinh(2D+2\epsilon)-2D-2\epsilon)}{\vol(M)}$.
If $\Sigma$ is a pseudo-minimal surface in $M'$, $\Sigma$
intersects at most $\kappa\left|\chi(\Sigma)\right|$ translates of
$\mathcal{D'}$ under the action of the group
$G=\pi_1(M)/\pi_1(M')$. From another viewpoint, for a given
translate of $\mathcal{D'}$ in $M'$, there exist at most $\kappa
\left|\chi(\Sigma)\right|$ copies of $\Sigma$ under the action of
$G$ which intersect it.

\end{lm}

\noindent \underline{Proof of lemma \ref{lm-intersectionD'}.}
\nopagebreak

Lemma \ref{lm-intersectionD'} is straightforward from inequality
(\ref{estimation_L}) of lemma \ref{lm-estimation_Dalpha&L}. The
embedded surface $\Sigma$ in $M'$ can be covered by at most
$a'\left|\chi(\Sigma)\right|$ embedded balls in $M'$ of radius
$\epsilon$. Therefore, this surface cannot intersect more than
$\lfloor \frac{\pi(\sinh(2D+2\epsilon)-2D-2\epsilon)}
{\vol(M)}a'\left|\chi(\Sigma)\right| \rfloor \leq
\frac{\pi(\sinh(2D+2\epsilon)-2D-2\epsilon)}
{\vol(M)}a'\left|\chi(\Sigma)\right|$ translates of $\mathcal{D'}$
in $M'$.\qed

~

Lemma \ref{lm-intersectionD'} applies to the pseudo-minimal
surface $\overline{F}$. Thus, this surface intersects at most
$\kappa\left|\chi(\overline{F})\right|\leq \kappa\,
\chi_-^c(\alpha')^2$ translates of $\mathcal{D'}$ in $M'$. Let $B$
be the subset of the corresponding elements of $G$.

Let also $C$ be the subset of $G$ corresponding to the translates
of $\mathcal{D'}$ that intersect $R'$. Still from lemma
\ref{lm-intersectionD'}, $\left|C\right|\leq
\kappa\left|\chi(R')\right|=\kappa \|\alpha'\|$.

The following claim and its proof are adapted from the proof of
\cite[Lemma 13]{Lac1}.

\begin{claim}

Set $\ell:=\sqrt[4]{117\kappa^2/4} $.

If $\ell\,\chi_-^c(\alpha') \leq \sqrt[4]{d}$, under the action of
$G$, there are at least $m'=9\chi_-^c(\alpha')/2$ translates of
$R'$ that are disjoint and do not intersect $\overline{F}$.

\end{claim}

\noindent \underline{Proof of claim.}
\nopagebreak

By contradiction, suppose that the claim is false. Then, for each
$m'$-uplet $(g_1R',\ldots,\\ g_{m'}R')$ of translates of $R'$, at
least two of them intersect, or at least one of them intersects
$\overline{F}$. There exist $j$ and $k$, with $1\leq j<k\leq m'$,
$c_1$ and $c_2\in C$ such that $g_jc_1=g_kc_2$, or there exist
$b\in B$, $c_1\in C$ and $s$ such that $g_sc_1=b$. In the first
case, $g_k^{-1}g_j\in CC^{-1}$, and in the second case, $g_s\in
BC^{-1}$. This means that $G^{m'}$ is the union of the sets
$q_{jk}^{-1}(CC^{-1})$ et $p_s^{-1}(BC^{-1})$, where for $1\leq
j<k \leq m'$ and $1 \leq s \leq m'$, $q_{jk}$ and $p_s$ are the
applications
\begin{eqnarray*}
q_{jk}\,:\,G^{m'} &\rightarrow & G\\
(g_1,\ldots,g_{m'}) &\mapsto & g_k^{-1}g_j\\
p_s\,:\,G^{m'} &\rightarrow & G\\
(g_1,\ldots,g_{m'}) &\mapsto & g_s.
\end{eqnarray*}
The cardinality of $q_{jk}^{-1}(CC^{-1})$ is
$\left|G\right|^{m'-1}\left|CC^{-1}\right|$, and the cardinality
of $p_s^{-1}(BC^{-1})$ is $\left|G\right|^{m'-1}
\left|BC^{-1}\right|$. Thus,

\begin{eqnarray*}
\left|G\right|^{m'}&\leq &\left(_{\,2}^{m'}\right)
\left|G\right|^{m'-1}\left|C\right|^2 +m'
\left|G\right|^{m'-1}\left|C\right|\left|B\right|\\
d^{m'}&\leq & \left(_{\,2}^{m'}\right) d^{m'-1} (\kappa
\|\alpha'\|)^2+m' d^{m'-1} \kappa \|\alpha'\| \kappa
\,\chi_-^c(\alpha')^2.
\end{eqnarray*}
As $\|\alpha'\|=\left|\chi(R')\right|\leq \left|\chi(S')\right| =
\chi_-^c(\alpha')$, one has
\begin{equation}\label{ineq-inegasymptotiqueLac}
d\leq \frac{\kappa^2}{2}m'(m'-1)\chi_-^c(\alpha')^2+\kappa^2 m'
\chi_-^c(\alpha')^3.
\end{equation}

As $m'=9\chi_-^c(\alpha')/2$, this leads to
\begin{eqnarray*}
d &\leq& \frac{9\kappa^2}{4}
\chi_-^c(\alpha')(\frac{9\chi_-^c(\alpha')}{2}-1)\chi_-^c(\alpha')^2+\frac{9\kappa^2}{2}
\chi_-^c(\alpha')^4\\
&\leq & \frac{117\kappa^2}{8} \chi_-^c(\alpha')^4 - \frac{9\kappa^2}{4}\chi_-^c(\alpha')^3\\
&\leq & \frac{117\kappa^2}{8} \chi_-^c(\alpha')^4.
\end{eqnarray*}

But as $\ell^4=117\kappa^2/4$ and $\ell^4\, \chi_-^c(\alpha')^4
\leq d$, one gets $d\leq d/2$, which provides a contradiction.
Therefore, the claim is true under those assumptions.\qed

~

From the claim, there exist at least $9\chi_-^c(\alpha')/2$
translates of $R'$ such that any two of them are disjoint, and
which do not intersect the surface $\overline{F}$ either. As each
of those $9\chi_-^c(\alpha')/2$ incompressible surfaces is in the
complement of $\overline{F}$, which is a disjoint union of
compression bodies, this surface is in fact parallel to a
component of $\overline{F}$. From corollary \ref{cor-estimates},
$\overline{F}$ has at most $3\chi_-^c(\alpha')/2$ components.
Therefore, there are at least three disjoint translates of $R'$
that are parallel to the same component of $\overline{F}$. Thus,
those three translates are parallel. If the surface $R'$ is
arbitrarily given an orientation,each of the translates of $R'$ is
oriented, and its orientation is given by the orientation of $R'$.
With those conventions, there are at least two of those parallel
translates whose orientations are coherent in the product region
they bound. Thus, there exists an incompressible surface $R''$ in
$M'$ and $h \in G$ an orientation preserving homeomorphism such
that $R''$ and $h(R'')$ are parallel and disjoint in $M'$. As
$R''$ is incompressible, Lemma 14 of \cite{Lac1} applies: the
cover $M'$ fibers over the circle, with fiber $R''$. But as $R''$
is a translate of the surface $R'$ under the action of $G$, if
$p\,:\,M'\rightarrow M$ is the covering map, the homology class of
$p^{-1}(R)$ is fibered. From a result of Gabai \cite[Lemme
2.4]{Ga1}, the homology class of $R$ is also fibered in $M$. As
$R$ is an embedded and incompressible surface (as it is
$\|\alpha\|$-minimizing), this means that the manifold $M$ fibers
over the circle, and that $R$ is a fiber.

There remains to show that if we do not a priori suppose that the
surface $R'$ is $\|\alpha'\|$-minimizing and such that
$h(M',\alpha',R')=h(M',\alpha')$, the surfaces $R'$ and $R$ are
still fibers. If $R''$ is a $\|\alpha'\|$-minimizing embedded
surface, such that $h(M',\alpha',R'')=h(M',\alpha')$, then the
argument above shows that $M'$ fibers over the circle and that
$R''$ is a fiber. But as the surface $R'$ is a component of the
preimage of $R$, it is incompressible and in the same homology
class as $R''$. Thus it is also a fiber. The argument above then
applies to show that $R$ is also a fiber. This ends the proof of
theorem \ref{thm-deccircregulier&fibration}.\qed

~

\noindent \underline{Proof of corollary
\ref{cor-deccircregulier&fibration-asymptotique}.}
\nopagebreak

The proof is straightforward from theorem
\ref{thm-deccircregulier&fibration}. If
$\lim_{i\div}\frac{\chi_-^c(\alpha_i)}{\sqrt[4]{d_i}}=0$, for $i$
large enough, $\ell\,\chi_-^c(\alpha_i) \leq \sqrt[4]{d_i}$, and
theorem \ref{thm-deccircregulier&fibration} applies.\qed

~

\noindent \underline{Proof of corollary
\ref{cor-revcycliquesduauxetcl}.}
\nopagebreak

As the cover $M_i\rightarrow M$ is the $i$-sheeted cyclic cover
associated to the class $\alpha$ and $\alpha_i=p_i^*(\alpha)$,
$\|\alpha_i\|=\|\alpha\|$. Thus,
$\chi_-^c(\alpha_i)=\|\alpha_i\|+2h(\alpha_i) =\|\alpha\|+2
h(\alpha_i)$. 
If there exists $i\geq i_0=\lceil (2\ell\|\alpha\|)^4\rceil$ such that
$\frac{h(\alpha_i)}{\sqrt[4]{i}}\leq \frac{1}{4\ell},$
then 
\begin{align*}
 \ell \chi_-^c(\alpha_i) &= \ell (\|\alpha\|+2
h(\alpha_i)) \leq \ell \| \alpha \| + \sqrt[4]{i}/2 \leq \sqrt[4]{i_0}/2 +
\sqrt[4]{i}/2 \leq \sqrt[4]{i}.
\end{align*}

Theorem \ref{thm-deccircregulier&fibration} then applies.\qed

\bibliographystyle{short}

\bibliography{biblio}

Claire \textsc{Renard},

\'Ecole Normale Sup\'erieure de Cachan,

Centre de Math\'ematiques et de Leurs Applications.

61 avenue du pr\'esident Wilson

F-94235 CACHAN CEDEX.

\emph{claire.renard@normalesup.org}

\end{document}